\documentclass[draft|showlink]{siamart1116}

\usepackage{amsfonts}
\usepackage{graphicx}
\usepackage{epstopdf}
\usepackage{algorithmic}
\ifpdf
 \DeclareGraphicsExtensions{.eps,.pdf,.png,.jpg}
\else
 \DeclareGraphicsExtensions{.eps}
\fi

\newtheorem{assumption}{Assumption}
\newtheorem{scheme}{Scheme}

\newcommand{\TheTitleABB}{generalized scheme for BSDEs and its error estimates} 
\newcommand{\TheTitle}{A generalized scheme for BSDEs based on derivative approximation and its error estimates} 
\newcommand{\TheAuthors}{Chol-Kyu Pak, Mun-Chol Kim and Hun O}

\def\R{\mathbb R} \def\P{\mathbb P} \def\F{\mathcal F} \def\d{\partial} \def\to{\rightarrow} 
\def\beq{\begin{equation}} \def\enq{\end{equation}}
\def\beseq{\begin{subequations}} \def\enseq{\end{subequations}}
\def\beqa{\begin{eqnarray}} \def\enqa{\end{eqnarray}}
\def\non{\nonumber}
\def\BeDef{\begin{definition}} \def\EnDef{\end{definition}}
\def\BeThe{\begin{theorem}} \def\EnThe{\end{theorem}}
\def\BeLem{\begin{lem}} \def\EnLem{\end{lem}}

\numberwithin{theorem}{section}
\numberwithin{equation}{section}
\numberwithin{table}{section}
\numberwithin{figure}{section}
\numberwithin{remark}{section}
\numberwithin{example}{section}

\headers{\TheTitleABB}{\TheAuthors}

\title{\TheTitle}

\author{
 Chol-Kyu Pak\thanks{Faculty of Mathematics, Kim Il Sung University, Pyongyang, Democratic People's Republic of Korea (\email{pck2016217@gmail.com}).}
\and
 Mun-Chol Kim\footnotemark[1]
\and
 Hun O\footnotemark[1]
}

\usepackage{amsopn}

\begin{document}

\maketitle

\begin{abstract}
In this paper we propose a generalized numerical scheme for backward stochastic differential equations(BSDEs). The scheme is based on approximation of derivatives via Lagrange interpolation. By changing the distribution of sample points used for interpolation, one can get various numerical schemes with different stability and convergence order. We present a condition for the distribution of sample points to guarantee the convergence of the scheme. \end{abstract}



\section{Introduction}\label{sec1}
Let $(\Omega,\F, \P)$ be a probability space, $T>0$ a finite time and $\{\F_t\}_{0\le t \le T}$ a filtration satisfying the usual conditions. Let $(\Omega,\F, \P,\{\F_t\}_{0\le t\le T})$ be a complete filtered probability space on which a standard $d$-dimensional Brownian motion $W_t=(W_t^1,W_t^2, \cdots ,W_t^d )^T$ is defined and $\F_0$ contain all the $\P$-null sets of $\F$. 
\par The standard form of backward stochastic differential equation (BSDE) is
\beq\label{bsde}
y_t=\xi+\int_t^T{f(s,y_s,z_s )ds}-\int_t^T{z_sdWs}, \quad t\in[0,T] 
\enq
where the generator $f:[0,T]\times\R^m\times\R^{m\times d}\to \R^m$ is $\{\F_t\}$-adapted for each $(y,z)$ and the terminal variable $\xi$ is a $\F_T$-measurable and square integrable random variable.
A process $(y_t,z_t ):[0,T]\times\Omega \rightarrow \R^m\times \R^{m\times d}$ is called an $L^2$-solution of the BSDE {\cref{bsde} if it satisfies the equation \cref{bsde} while it is ${\{\F_t\}}$-adapted and square integrable .
\par In 1990, Pardoux and Peng proved in \cite{Peng90} the existence and uniqueness of the solution of general nonlinear BSDEs and afterwards there has been active works in this field. BSDEs and their extensions are very important in mathematical finance(\cite{Karoui97}), stochastic control(\cite{Bismut73, Hamadene95}) and partial differential equations (\cite{Peng91, Pardoux98}).
\par In this paper we assume that the terminal condition is a function of $W_T$, i.e., $\xi=\varphi(W_T)$ and the BSDE \cref{bsde} has a unique solution $(y_t,z_t )$. 
It was shown in \cite{Peng91} that the solution $(y_t,z_t)$ of \cref{bsde} can be represented as
\beq\label{eq12}
y_t=u(t,W_t ), z_t=\nabla_x u(t,W_t ), \; \forall t\in[0,T)
\enq
where $u(t,x)$ is the solution of the parabolic partial differential equation
\beq\label{pde}
\frac{\d u}{\d t}+\frac{1}{2} \sum_{i=1}^d \frac{\d^2 u}{\d x_i^2}+f(t,u,\nabla_x u)=0
\enq
with the terminal condition $u(T,x)=\varphi (x)$, and $\nabla_x u$ is the gradient of $u$ with respect to the spatial variable $x$. The smoothness of $u$ depends on $f$ and $\varphi$.
\par Although BSDEs and their extensions such as FBSDEs have very important applications in many fields, it is well known that it is difficult to obtain the analytic solutions except some special cases and there have been many works on numerical methods. (\cite{Bou, Gobet, Prot02, Prot94, Tre06, Zhao06, Ruijter, Zhao14, Zhao17}) 

By making use of the relations between the BSDEs and the partial differential equations, a four step scheme was proposed in \cite{Prot94, Prot02} and a layer method for PDEs was applied to the BSDEs to the random time horizon in \cite{Tre06}. On the other hand several kinds of probabilitistic approaches have also been developed and they generally need approximation of conditional expectations. Bouchard and Touzi\cite{Bou} developed a backward simulation scheme and Gobet\cite{Gobet} proposed a method based on regression on function bases. Ruijter and Oosterlee\cite{Ruijter} developed a Fourier cosine method and Zhao and his co-authors have developed several schemes for various types of BSDEs and FBSDEs using Gauss-Hermite quadrature rule for approximation of conditional expectations. \cite{Zhao06, Zhao10, Zhao14, Zhao17}

\par Among all the previous works, our work is closely related to \cite{Zhao14} where a kind of multi-step scheme for FBSDEs was proposed and Gauss-Hermite quadrature rule was used for approximation of conditional expectations. They derived a new kind of reference scheme which was much simpler than the previous ones. They demonstrated experimentally that their scheme is of up to 6th order and stated that their scheme does not seem to converge for higher order because the coefficients fail to satisfy root condition.
\par Our scheme is based on the same reference equation with \cite{Zhao14} but use general form of derivative approximation. We got an inspiration of nonequidistant difference scheme while trying to avoid spatial interpolation for the approximation of conditional expectation using Gauss-Hermite quadrature rule. In the previous works spatial interpolation has been inevitable because it was impossible to set up a fully nested time-space grid and it contributed much to complexity. But through some investigation, we showed that if we derive a scheme using quadratically distant derivative approximation we can set up a fully nested time-space grid for second or third order and avoid spatial interpolation. (\cite{PCK18})
\par In this paper we generalized a numerical scheme based on derivative approximation which varies according to the distribution of sample points used for derivative approximation. We present a condition for the "parameters" under which the resulting scheme converges. The result on error estimates supports the previous works \cite{Zhao14, PCK18}.
\par We note that our scheme is limited to the case of Brownian setting, Markovian terminal condition and deterministic time horizon. One can refer to \cite{Tre06} for numerical solution of random terminal time problems. 

\par The rest of the paper is organized as follows. In \Cref{sec2}, we introduce a generalized discrete scheme for BSDEs. In \Cref{sec3}, we present the error estimates of the scheme for semidiscrete version for certain types of generator function. In \Cref{sec4}, we discuss some specific resulting schemes. Finally some conclusions are given in \Cref{sec5}.

\section{A generalized numerical scheme for BSDEs based on derivative approximation}\label{sec2}
\subsection{The reference equations}\label{subsec21}
We firstly note that the major content of this subsection is from \cite{Zhao14} where decoupled FBSDEs were discussed.
\par Let us consider the following backward stochastic differential equation:
\beq\label{eq21} y_t=\varphi(W_T)+\int_t^T{f(s,y_s,z_s)ds}-\int_t^T{z_sdW_s}, \quad t\in [0,T] \enq
where the generator $f:[0,T]\times\Omega\times\R^m\times\R^{m\times d}\to\R^m$ is a stochastic process that is $\{\F_t\}$-adapted for all $(y,z)$ and $\varphi:\R^d\to\R^m$ is a measurable function. 
\par Throughout the paper, we make assumptions on the parameters as follows.
\begin{assumption}\label{ass1}
The generator function $f$ and $\varphi$ are smooth enough and their derivatives as well as themselves are all bounded.
\end{assumption}
\par Let $0=t_0<\cdots<t_N=T$ be an equidistant partition of $[0,T]$ and $t_{n+1}-t_n=h=T/N$.
As in \cite{Zhao06,Zhao10}, we define $\F_s^{t,\mathbf{x}}(t\leq s\leq T)$ to be be a $\sigma$-field generated by the Brownian motion $\{\mathbf{x}+W_r-W_t,t\leq r\leq s\}$ starting from the time-space point $(t,\mathbf{x})$ and $E_t^\mathbf{x}[\cdot]:=E[\cdot|\F_t^{t,\mathbf{x}}]$.
\par Taking the conditional expectation $E_{t_n}^\mathbf{x}[\cdot]$ on the both sides of \cref{eq21}, we get
\beq\label{eq22} E_{t_n}^\mathbf{x}[y_t]=E_{t_n}^\mathbf{x}[\varphi(W_T)]+\int_t^T{E_{t_n}^\mathbf{x}[f(s,y_s,z_s)]ds}\enq
By differentiating the both sides of \cref{eq22} with respect to $t$, we obtain the following ODE.(see \cite{Zhao14} for differentiability)
\beq\label{eq23} \frac{dE_{t_n}^\mathbf{x}[y_t]}{dt}=-E_{t_n}^\mathbf{x}[f(s,y_s,z_s)]\enq
Next, for $t\in[t_n,T]$ we have
\beq\label{eq24} y_{t_n}=y_t+\int_t^T{f(s,y_s,z_s)ds}-\int_{t_n}^t{z_sdW_s}.\enq
Multiplying $\Delta W_{t_n,t}^\mathbf{T}:=W_t^\mathbf{T}-W_{t_n}^\mathbf{T}$ to the both sides of \cref{eq24} and taking $E_{t_n}^\mathbf{x}[\cdot]$, we obtain
\beq\label{eq25} 0=E_{t_n}^\mathbf{x}[y_t\Delta W_{t_n,t}^\mathbf{T}]+\int_{t_n}^t{E_{t_n}^\mathbf{x}[f(s,y_s,z_s)\Delta W_{t_n,t}^\mathbf{T}]ds}-\int_{t_n}^t{E_{t_n}^\mathbf{x}[z_s]ds}\enq
where $(\cdot)^\mathbf{T}$ means transpose of $(\cdot)$.
\par Taking derivatives of the both sides of \cref{eq25} with respect to $t$, we obtain the following ODE.
\beq\label{eq26} E_{t_n}^\mathbf{x}[z_t]=\frac{dE_{t_n}^\mathbf{x}[y_t\Delta W_{t_n,t}^\mathbf{T}]}{dt}+E_{t_n}^\mathbf{x}[f(t,y_t,z_t)\Delta W_{t_n,t}^\mathbf{T}]\enq
The two ODEs \cref{eq23} and \cref{eq26} are called the reference equations for the BSDE \cref{eq21}.
From these reference equations we have
\beq\label{eq27} \frac{dE_{t_n}^\mathbf{x}[y_t]}{dt}|_{t=t_n}=-f(t_n,y_{t_n},z_{t_n})\enq
\beq\label{eq28} z_{t_n}=\frac{dE_{t_n}^\mathbf{x}[y_t\Delta W_{t_n,t}^\mathbf{T}]}{dt}|_{t=t_n}.\enq

\subsection{Derivative approximation and semi-discrete scheme for BSDE}
Now by approximating the derivatives in \cref{eq27} and \cref{eq28} we get discrete scheme for BSDE \cref{eq21}. Among many methods for numerical differentiation, we will use the one based on Lagrange interpolation. In \cite{Zhao14} the approximation formula was derived by solving linear system obtained by Taylor's expansion and it coincides with the one from Lagrange interpolation based on equidistant points. 
\par Let $u(t):\R\to\R$ be $k+1$ times differentiable and $a_0=0<a_1<a_2<\cdots<a_k<N$ be $k+1$ integers, then the Lagrange interpolation polynomial based on values $\{u_0,u_1,\cdots,u_k\}$ on sample points $\{t_0,t_0+a_1h,\cdots,t_0+a_kh\}$ can be written as
\beq\label{eq29} L(t)=\sum_{i=0}^k{\frac{\prod_{i\neq j}{(t-t_0-a_jh)}}{\prod_{i\neq j}{(a_i-a_j)}}h^{-k}u_i}\enq
and the deviation is given by
\beq\label{eq210} L(t)-u(t)=\frac{f^{(k+1)}(\xi)\prod_{i=0}^n{(t-t_0-a_ih)}}{(k+1)!}, \quad t_0\leq \xi \leq t_0+a_kh.\enq
Hereafter we will call $a_1,\cdots,a_k$ the "parameters of the scheme".
By differentiating \cref{eq29}, we get
\beq\label{eq211} \frac{dL}{dt}(t)=\sum_{i=0}^k{\frac{h^{-k}u_i}{\prod_{j\neq i}{(a_i-a_j)}}\sum_{j\neq i}{\prod_{l\neq i, l\neq j}{(t-t_0-a_lh)}}}\enq
and furthermore,
\beq\label{eq212}
\frac{du}{dt}(t_0)=\frac{\sum_{i=0}^k{\gamma_i^ku_i}}{h}+\frac{(-1)^ku^{(k+1)}(\xi)}{k}h^k
\enq
where 
\beq\label{eq213}
\gamma_0^k=-\sum_{j\neq 0}{\frac{1}{a_j}}, \gamma_i^k=\frac{(-1)^{k-1}a_1a_2\cdots a_k}{a_i\prod_{i\neq j}{(a_i-a_j)}}, (1\leq i \leq k).
\enq
So under the assumption that $u^{(k+1)}$ is bounded we have
\beq\label{eq214}
\frac{du}{dt}(t_0)=\frac{\sum_{i=0}^k{\gamma_i^ku_i}}{h}+O(h^k)
\enq
\\
About the coefficients, the following proposition holds.
\bigskip
\begin{proposition}\label{pro21}
The coefficients $\gamma_i^k$ satisfy 
$\sum_{i=0}^k{\gamma_i^k}=0$
\end{proposition}
\bigskip
Now we apply this approximation formula to the reference equations \cref{eq27} and \cref{eq28} to get 
\beq\label{eq215}
-\gamma_0^ky_{t_n}=\sum_{i=1}^k{\gamma_i^kE_{t_n}^x[y_{t_{n+a_i}}]}+hf(t_n,y_{t_n})+hR_{y,n}^k
\enq
\beq\label{eq216}
hz_{t_n}=\sum_{i=1}^k{\gamma_i^kE_{t_n}^x[y_{t_{n+a_i}}\Delta W_{n,a_i}^\mathbf{T}]}+hR_{z,n}^k
\enq
where $R_{y,n}^k, R_{z,n}^k$ are the truncation errors as follows. 
\[
R_{y,n}^k={\frac{dE_{t_n}^x[y_t]}{dt}\vline}_{t=t_n} - \sum_{i=0}^k{\frac{\gamma_i^k}{h}E_{t_n}^x[y_{t_{n+a_i}}]}\]
\[
R_{z,n}^k={\frac{dE_{t_n}^x[y_t\Delta W_{t_n,t}]}{dt}\vline}_{t=t_n} - \sum_{i=0}^k{\frac{\gamma_i^k}{h}E_{t_n}^x[y_{t_{n+a_i}}\Delta W_{t_n,t}^\mathbf{T}]}.
\]
Note that we denote $W_{t_{n+a_i}}-W_{t_n}$ by $\Delta W_{n,a_i}$ for the simplicity in notion.
From \cref{eq214} one can easily prove the following lemma.
\bigskip
\begin{lemma}\label{lem21}
Under the assumption \cref{ass1}, the truncation erros satisfy 
\[ |R_{y,n}^k|\leq Ch^k, |R_{z,n}^k|\leq Ch^k \]
where $C$  is a constant that depends only on bounds of $f,\varphi$ and their derivatives.
\end{lemma}
\bigskip
Now we present a semi-discrete scheme for BSDE \cref{eq21}.
\bigskip
\begin{scheme}\label{scheme1}
Assume that $y^{N-i},z^{N-i}(i=0,\cdots,a_k-1)$ are known. For $n=N-a_k,\cdots,0$ solve $y^n=y^n(\mathbf{x}),z^n=z^n(\mathbf{x})$ at time-space point $(t_n,\mathbf{x})$ by
\beq\label{eq217} -\gamma_0^ky^n=\sum_{j=1}^k{\gamma_j^kE_{t_n}^\mathbf{x}[y^{n+a_j}]}+h\cdot f(t_n,y^n,z^n)
\enq

\beq\label{eq218} z^n=\frac{\sum_{j=1}^k{\gamma_j^kE_{t_n}^\mathbf{x}[y^{n+a_j}\Delta W_{n,a_j}^\mathbf{T}]}}{h}\enq
\end{scheme}
\bigskip
\par This is called a semi-discrete scheme because it is discretized only in time domain and one needs a spatial grid to get the fully discrete scheme. We will discuss about the full discrete scheme later in \Cref{sec4}.
\par Note that the scheme is explicit for $z$ but implicit for $y$. At every step we calculate $z$ component first and then get $y$ through Picard's iteration.
\par Note that one can make various schemes by selecting $a_1,\cdots,a_k$ variously. For the case where $a_j=j (j=1,\cdots, k)$ the scheme is exactly the same with the one discussed in \cite{Zhao14}. If we select $a_j=j^2 (j=1,\cdots, k)$ we get a one based on quadratically distant sample points which is considered 
\section{Error estimates}\label{sec3}
In this section we prove the error estimates of the semidiscrete \Cref{scheme1} for the case where the generator is independent on control variable $z$. We assume $m=d=1$ for the sake of simplicity.
\beq\label{eq31}
y_t=\varphi(W_T)+\int_t^T{f(s,y_s)ds}-\int_t^T{z_sdW_s}, t\in[0,T]
\enq
\par Let $y_t, z_t$ be the solutions of the BSDE \cref{eq31} and $y^n,z^n$ be the solutions of the \Cref{scheme1}. We define error terms like $e_y^n:=y_{t_n}-y^n, e_z^n:=z_{t_n}-z^n$.
\par We also define $c_j, \alpha_j (j=0,1,\cdots,a_k)$ for the statement of the estimates as follows.
\beq\label{eq32}
 c_j :=
   \left\{ \begin{array}{l}
       \gamma_l^k,  \exists 0\leq l\leq k, j=a_l \\ 
	0, otherwise \end{array} \right. ,  \alpha_j := \sum_{l=j}^{a_k}{c_l}
\enq
\par The following theorem gives the error estimates for $y$.
\bigskip
\begin{theorem}\label{the31}
Suppose that \Cref{ass1} holds and the initial values satisfy $\max\limits_{N-a_k+1\leq n\leq N}{E[|e_y^n|]}=O(h^k)$. Furthermore suppose that the numbers $\alpha_1, \cdots, \alpha_{a_k}$ that are defined as \cref{eq32} satisfy the following inequality.
\beq\label{eq33}\frac{\sum_{j=2}^{a_k}{|\alpha_j|}}{|\alpha_1|}<1 \enq
Then for sufficiently small time step $\Delta t$ we have
\beq\label{eq34} \max\limits_{0\leq n\leq N-a_k}{E[|e_y^n|}\leq Ch^k \enq
where $C>0$ is a constant depending only on $\{\gamma_j\}, T$, the upper bounds of functions $\varphi$ and $f$, and their derivatives.
\end{theorem}
\begin{proof}
\medskip
We prove the theorem in two steps.
\par \textbf{STEP 1.} We assume that the time horizon $T$ is small, say more accurately 
\beq\label{eq35}
T<\frac{|\alpha_1|-\sum_{i=2}^{a_k}|\alpha_i|}{2L}
\enq
where $L$ is a Lipschitz constant of the generator function $f$.\\
From the reference equation \cref{eq27} and the scheme, we can write
\[
-\gamma_0^ke_y^n=\sum_{j=1}^k{\gamma_j^kE_{t_n}^x[e_y^{n+a_j}]}+h\Delta f_n+O(h^{k+1})\]
where \[\Delta f_n = f(t_n,y_{t_n})-f(t_n,y^n).\]
From the \Cref{pro21}, $\sum_{j=0}^k{\gamma_j^k}=0$ and we have
\[
\sum_{j=1}^k{\gamma_j^ke_y^n}=\sum_{j=1}^k{\gamma_j^kE_{t_n}^x[e_y^{n+a_j}]}+h\Delta f_n+O(h^{k+1})
\]
which leads to 
\[
\sum_{j=1}^k{\gamma_j^k(e_y^n-E_{t_n}^x[e_y^{n+a_j}])}=h\Delta f_n+O(h^{k+1}).
\]
Using $c_j$ we can rewrite the above equation as follows.
\beq\label{eq36}
\sum_{j=1}^{a_k}{c_j(e_y^n-E_{t_n}^x[e_y^{n+j}])}=h\Delta f_n+O(h^{k+1}).
\enq
Similarly we can get the following equations.
\beseq
\begin{align}
\sum_{j=1}^{a_k}{c_j(E_{t_n}^x[e_y^{n+1}]-E_{t_n}^x[e_y^{n+1+j}])}&=h\Delta f_{n+1}+O(h^{k+1})\non\\
&\vdots\non\\
\sum_{j=1}^{a_k}{c_j(E_{t_n}^x[e_y^{N-a_k}]-E_{t_n}^x[e_y^{N-a_k+j}])}&=h\Delta f_{N-a_k}+O(h^{k+1})\non
\end{align}
\enseq
Summing up all the above equations including \cref{eq36} we have
\beseq
\begin{align}
\alpha_1 e_y^n + &\sum_{j=2}^{a_k}{E_{t_n}^x[e_y^{n+j-1}]}-\sum_{j=1}^{a_k}{\alpha_iE_{t_n}^x[e_y^{N+i-a_k}]}=\non\\
&h\Delta f_n + h\sum_{j=n+1}^{N-a_k}{E_{t_n}^x[\Delta f_j]} + (N-a_k-n+1)O(h^{k+1})\non
\end{align}
\enseq
where $\alpha_1, \cdots, \alpha_{a_k}$ that are defined as \cref{eq32}.
\par Taking expectation to the both sides of the above equation, we have
\[
|\alpha_1|E[|e_y^n|]\leq\sum_{i=2}^{a_k}{|\alpha_i|E[|e_y^{n+i-1}|]}+\sum_{i=1}^k{|\alpha_i|E[|e_y^{N+i-a_k}|]}+h\sum_{j=n}^{N-a_k}{E[|\Delta f_j|]} + NO(h^{k+1}).
\]
Then from the Lipschitz continuity of the generator function and the assumption on the initial approximation, there exists a constant $C$ such that
\[
|\alpha_1|E[|e_y^n|]\leq\sum_{i=2}^{a_k}{|\alpha_i|E[|e_y^{n+i-1}|]}+Lh\sum_{j=n}^{N-a_k}{E[|e_y^j|]} + Ch^k.
\]
We used the fact that $N=T/h$.
\par Rearranging the terms leads to
\[
(1-\frac{Lh}{|\alpha_1|})E[|e_y^n|]\leq \sum_{i=2}^{a_k}{\frac{|\alpha_i|}{|\alpha_1|}E[|e_y^{n+i-1}|]}+\frac{Lh}{|\alpha_1|}\sum_{j=n+1}^{N-a_k}{E[|e_y^j|]} + Ch^k.
\]
Let $p:=\arg\max\limits_{0\leq n\leq N-a_k}{E|e_y^n|}$ then
\beseq
\begin{align}
(1-\frac{Lh}{|\alpha_1|})E[|e_y^p|]&\leq\sum_{i=2}^{a_k}{\frac{|\alpha_i|}{|\alpha_1|}E[|e_y^{p+i-1}|]}+\frac{Lh}{|\alpha_1|}\sum_{j=p+1}^{N-a_k}{E[|e_y^j|]} + Ch^k\non\\
&\leq E[|e_y^p|]\sum_{i=2}^{a_k}{\frac{|\alpha_i|}{|\alpha_1|}}+\frac{Lh}{|\alpha_1|}E[|e_y^p|](N-a_k-p)+Ch^k\non\\
&\leq E[|e_y^p|]\left(\sum_{i=2}^{a_k}{\frac{|\alpha_i|}{|\alpha_1|}}+\frac{LT}{|\alpha_1|}\right)+Ch^k\non
\end{align}
\enseq
So we get
\[
\left(1-\sum_{i=2}^{a_k}{\frac{|\alpha_i|}{|\alpha_1|}}-\frac{LT}{|\alpha_1|}(1+\frac{1}{N})\right)E[|e_y^p|]\leq Ch^k.\]
From \cref{eq33} and \cref{eq35} we finally get
\[E[|e_y^p|]\leq Ch^k\]
where $C$ is a generic constant which only depends on $\{\gamma_j\}, T$, the upper bounds of functions $\varphi$ and $f$, and their derivatives.

\par  \textbf{STEP 2.} Now let us consider the case $T$ is big. In this case we can divide the time interval $[0,T]$ into several subintervals that satisfy \cref{eq35}. Then we can get the result through typical procedure. 
\end{proof}
\bigskip
Now we discuss the error for $z$. We introduce the variational equation of \cref{eq31} as follows.
\beq\label{eq310} \nabla y_t=\varphi_x(W_T)+\int_t^T{f_y(s,y_s)\nabla y_sds}-\int_t^T{\nabla z_sdW_s}\enq
where  $\varphi_x,f_y$ are the partial derivatives of $\varphi,f$  with respect to  $x,y$ and $\nabla y_s,\nabla z_s$ are the variations of $y_s,z_s$ with respect to spatial variable  $x$.(See \cite{Zhao10}) From the relationship \cref{eq12}, $z_t=\nabla y_t.$\\
For the variational equation we can derive the following equation similarly as in \cref{eq215} and \cref{eq216}.
\beq\label{eq311}
-\gamma_0^k\nabla y_{t_n} = \sum_{i=0}^k{\gamma_i^k}E_{t_n}^x[\nabla y_{t_{n+a_i}}]+hf_y(t_n,y_{t_n})\nabla y_{t_n}+hR_{\nabla y,n}^k
\enq
where
\[R_{\nabla y,n}^k={\frac{dE_{t_n}^x[\nabla y_t]}{dt}\vline}_{t=t_n}-\sum_{i=0}^k{\frac{\gamma_i^k}{h}E_{t_n}^x[\nabla y_{t_{n+a_i}}]}.\]
And from \cref{eq217}, we also have the following equation for the variation of $y^n$.
\beq\label{eq312}
-\gamma_0^k\nabla y_{n} = \sum_{i=0}^k{\gamma_i^k}E_{t_n}^x[\nabla y_{{n+a_i}}]+hf_y(t_n,y_{n})\nabla y_{n}.
\enq
The following lemma is about the bound of $y_t$ and $\nabla y_t$. (See Lemma 4 in \cite{Zhao10} and references therein for details.)
\bigskip
\begin{lemma}\label{lem32}
Suppose that \Cref{ass1} holds, and let $y_t$ and $\nabla y_t$ be the solutions of \cref{eq31} and \cref{eq310}, respectively. Then it holds that
\beq\label{eq313}
\sup\limits_{t\in[0,T]}{E[|y_t|]}+\sup\limits_{t\in[0,T]}{E[|\nabla y_t|]}\leq C ,
\enq
where $C>0$ is a generic constant depending only on $T$, the upper bounds of functions $\varphi$ and $f$, and their derivatives.
\end{lemma}
\bigskip
We define $e_{\nabla y}^n:=\nabla y_{t_n} - y^n$. Then through the similar procedure as in the \Cref{the31} we can prove the following lemma about $e_{\nabla y}^n$.
\bigskip
\begin{lemma}\label{lem33}
Suppose that \Cref{ass1} holds and the initial approximations satisfy $\max\limits_{N-a_k+1\leq n\leq N}{E[|e_\nabla y^n|]}=O(h^k)$. 
\\Furthermore suppose that the following inequality holds.
\beq\label{eq314}\frac{\sum_{j=2}^{a_k}{|\alpha_j|}}{|\alpha_1|}<1 \enq
Then for sufficiently small time step $\Delta t$ we have
\beq\label{eq315} \max\limits_{0\leq n\leq N-a_k}{E[|e_\nabla y^n|}\leq Ch^k \enq
where $C>0$ is a constant depending only on $\{\gamma_j\}, T$, the upper bounds of functions $\varphi$ and $f$, and their derivatives.
\end{lemma}
\bigskip
Now we present a theorem which gives the error estimates for $z$.
\bigskip
\begin{theorem}\label{the32}
Suppose that \Cref{ass1} holds and the initial approximations satisfy $\max\limits_{N-a_k+1\leq n\leq N}{E[|e_y^n|]}=O(h^k)$ and $\max\limits_{N-a_k+1\leq n\leq N}{E[|e_z^n|]}=O(h^k)$. Furthermore suppose that the numbers $\alpha_1, \cdots, \alpha_{a_k}$ that are defined as \cref{eq32} satisfy the following inequality.
\beq\label{eq316}\frac{\sum_{j=2}^{a_k}{|\alpha_j|}}{|\alpha_1|}<1 \enq
Then for sufficiently small time step $\Delta t$ we have
\beq\label{eq317} \max\limits_{0\leq n\leq N-a_k}{E[|e_z^n|}\leq Ch^k \enq
where $C>0$ is a constant depending only on $\{\gamma_j\}, T$, the upper bounds of functions $\varphi$ and $f$, and their derivatives.
\end{theorem}
\begin{proof}
From \cref{eq216} and \cref{eq218} we have
\beq\label{eq318}
he_z^n=\sum_{j=1}^{a_k}{\gamma_j^k E_{t_n}^x[e_y^{n+a_j}\Delta W_{n,a_j}]}+O(h^{k+1}).
\enq
We can derive
\beseq
\begin{align}
E_{t_n}^x[e_y^{n+a_j}\Delta W_{n,j}]&=E[e_y^{n+a_j}(x+\Delta W_{n,a_j}) \Delta W_{n,a_j}]=\non\\ &=\int\limits_R{\frac{1}{\sqrt{2\pi a_j h}}exp(-\frac{y^2}{2a_jh})dy}=\non\\&=a_jhE[e_{\nabla y}^{n+a_j}(x+\Delta W_{n,a_j})\Delta W_{n,a_j}]=\non\\&=a_jhE_{t_n}^x[e_{\nabla y}^{n+a_j}]\non
\end{align}
\enseq
By simple substitution we have
\[
e_z^n=\sum_{j=1}^{a_k}{\gamma_j^ka_jE_{t_n}^x[e_{\nabla y}^{n+a_j}]}+O(h^k)
\]
and this leads to 
\[
E[|e_z^n|]\leq \sum_{j=1}^{a_k}{|\gamma_j^k|a_jE[|e_{\nabla y}^{n+a_j}|]}+Ch^k.
\]
where $C>0$ is a constant.
So from \Cref{lem33} we obtain 
\[
\max\limits_{N-a_k+1\leq n\leq N}{E[|e_z^n|]}\leq Ch^k,
\]
where $C>0$ is a constant depending only on $\{\gamma_j\}, T$, the upper bounds of functions $\varphi$ and $f$, and their derivatives.
\end{proof}
\bigskip
\section{Discussions}\label{sec4}
In \cite{Zhao10, Zhao14} the authors stated that the stability of \Cref{scheme1} is closely related to its corresponding characteristics polynomial defined as follows.
\beq\label{eq41}
P_\gamma^k(\lambda):=\sum_{i=0}^k{\gamma_i^k \lambda^{a_k-a_i}}
\enq
The roots $\{\lambda_i^k\}_{i=1}^k$ of this polynomial are said to satisfy the root conditions if $|\lambda_i^k|\leq 1$ and $|\lambda_i^k|=1\Rightarrow P_\gamma^k(\lambda_i^k)'\neq 0$.
In \cite{Zhao14}, they showed that the roots of $P_\gamma^k(\lambda)$ satisfy the root conditions for only $1\leq k\leq 6$ when the sample points are equally distant, i.e. $a_i=i$. In \cite{PCK18} they proposed a scheme based on quadratically distant sample points(i.e. $a_i=i^2$) and it was effective in decrease of the computational complexity.
\par Now \cref{the31} and \cref{the32} show that the convergence of the \Cref{scheme1} is also decided by the parameters $a_1,\cdots, a_k$. Let us test if the parameters satisfy the condition \cref{eq33} for two kinds of selection of parameters. The first one will be the case where $a_i=i$ and the second one will be the case where $a_i=i^2$.
In \Cref{tbl41} we present two kinds of test values for several schemes changing the value of $k$ from 2 to 7.
The first value is the maximum absolute values of the roots except 1 and it represents the stability of the scheme.
The second value is ${\sum_{i=2}^{a_k}{|\alpha_i|}}/{|\alpha_1|}$ which represents the convergence of the scheme. 
\begin{table}[tbhp]\label{tbl41}
\caption{Comparison of two schemes in stability and convergence}
\centering
\begin{small}
\begin{tabular}{|c|c|c|c|c|c|c|c|c|c|c|c|c|} \hline
Param&\multicolumn{2}{|c|}{$k=2$}&\multicolumn{2}{|c|}{$k=3$}&\multicolumn{2}{|c|}{$k=4$}&\multicolumn{2}{|c|}{$k=5$}&\multicolumn{2}{|c|}{$k=6$}&\multicolumn{2}{|c|}{$k=7$}\\ \hline
$a_i=i$ &  0.33 & 0.33 & 0.81 & 0.42& 1.56 & 0.56& 2.73 & 0.70& 4.65 & 0.86&7.87&1.02\\ \hline
$a_i=i^2$ &  0.06 & 0.48 & 0.20 & 0.63& 0.37 & 0.73& 0.44 & 0.80& 0.50 & 0.83& 0.57 & 0.85\\ \hline
\end{tabular}
\end{small}
\end{table}
As one can see the equidistant scheme ($a_i=i$) fails to satisfy the condition for convergence from $k=4$ while the nonequidistant scheme satisfy all conditions for all test cases.
This shows that although one would need more initial approximations to use nonequidistant difference scheme, it could have a better convergence and stability.
\section{Conclusions}\label{sec5}
In this paper we proposed a generalized numerical scheme for BSDEs. The scheme is based on approximation of derivatives via Lagrange interpolation on generally distributed sample points. By changing the distribution of sample points used for interpolation, one can get various numerical schemes with different stability and convergence order. We presented a condition for the distribution of sample points to guarantee the convergence of the scheme. Through discussion we showed that a scheme based on quadratically distant sample points is both stable and convergent.

\bibliographystyle{siamplain}

\end{document}